\documentclass[12pt]{amsart}

\usepackage[all,ps,cmtip,dvips]{xy}
\usepackage[dvips]{graphicx}
\usepackage{psfrag}
\usepackage{amssymb}
\usepackage{amsfonts}
\usepackage[mathcal]{eucal}

\newtheorem{theo}{Theorem}[section]
\newtheorem{prop}[theo]{Proposition}
\newtheorem{lemm}[theo]{Lemma}

\newtheorem{exam}[theo]{Example}
\newtheorem{rema}[theo]{Remark}
\newtheorem{coro}[theo]{Corollary}

\newcommand{\Sun}{{\mathbf{SU}(n)}}

\newcommand{\ba}[1]{\begin{array}{#1}}
\newcommand{\ea}{\end{array}}
\newcommand{\lra}{\longrightarrow}
\newcommand{\lmt}{\longmapsto}

\renewenvironment{proof}{
\mbox{}\\
\textit{Proof. }}
{\mbox{}\hspace{\stretch{1}}$\Box$ \\}

\begin{document}
\title{Kirwan map and moduli space of flat connections}
\author{S\'ebastien Racani\`ere}
\address{Mathematics Department\\
South Kensington Campus\\
Imperial College London\\
SW7 2AZ\\
UK}
\email{s.racaniere@ic.ac.uk}
\urladdr{www.ma.ic.ac.uk/\~{}racani}
\maketitle

\let\ldots\dots
\def\F{\mathbf{F}}
\def\FF{\mathbb{F}}
\def\NN{\mathbb{N}}
\def\ZZ{\mathbb{Z}}
\def\QQ{\mathbb{Q}}
\def\RR{\mathbb{R}}
\def\CC{\mathbb{C}}
\def\HH{\mathbb{H}}
\def\TT{\mathbf{T}}
\def\M{\mathcal{M}}
\def\Fr{\rm{ Frac}}
\def\In{\rm{ Int}}
\def\Im{{\rm Im}}
\def\Re{{\rm Re}}
\def\id{{\rm id}}
\def\Id{{\rm Id}}
\def\Ad{{\rm Ad}}
\def\ad{{\rm ad}}
\def\det{{\rm det}}
\def\tr{{\rm tr}}
\def\ov{\overline}
\def\deg{{\rm deg}}
\def\Hom{{\rm Hom}}
\def\End{{\rm End}}
\def\ev{{\rm ev}}
\def\Sta{{\rm Stab}}
\def\stab{{\rm stab}}
\def\Orb{{\rm Orb}}
\def\P{{\rmP}}
\def\Cp{{\rm\CC P}}
\def\d{{\rm d}}
\def\T{{\rm T}}
\def\Exp{{\rm Exp}}
\def\grad{{\rm grad}}
\def\infinity{\infty}
\def\equ{equivariant}
\def\coh{cohomology}
\def\fnck{{F_n(\CC^k)}}
\def\fncinf{{F_n(\CC^\infinity)}}
\def\gnck{{G_n(\CC^k)}}
\def\gncinf{{G_n(\CC^\infty)}}
\def\homotopic{{\sim_h}}
\def\vol{\mathcal{V}}
\def\NK{{\rm NK}}
\def\NKb{{\rm\overline{NK}}}
\def\NG{{\rm NG}}
\def\NGb{{\rm\overline{NG}}}
\def\del{{\partial}}

\def\BTW{ Bott, Tolman and Weitsman }
\def\ABBR{ Atiyah-Bott -- Biswas-Raghavendra }

\begin{abstract}
If $K$ is a compact Lie group and $g\geq 2$ an integer, the space $K^{2g}$ is endowed with the structure
of a Hamiltonian space with a Lie group valued moment map $\Phi$. Let $\beta$ be in the centre of $K$. The
reduction $\Phi^{-1}(\beta)/K$ is homeomorphic to a moduli space of flat connections.
When $K$ is simply connected, a direct consequence of a recent paper of Bott,
Tolman and Weitsman is to give a set of generators for the  $K$-equivariant cohomology of $\Phi^{-1}(\beta)$.

Another method to construct classes in $H^*_K(\Phi^{-1}(\beta))$ is by using the so called universal bundle. When
the group is $\Sun$ and $\beta$ is a generator of the centre, these last classes are known to also generate the
equivariant cohomology of $\Phi^{-1}(\beta)$. The aim of this paper is to compare the classes constructed
using the result of Bott, Tolman and Weitsman and the ones using the universal bundle.

In particular, I prove that the set of cohomology classes coming from the universal bundle is indeed
a set of multiplicative generators for the cohomology of the moduli space. With
$K=\Sun$, this is a new proof of the well-known construction of generators for the cohomology
of the moduli space of semi-stable vector bundles with fixed determinant.
\end{abstract}


\section*{Index of notation}

\begin{tabular}{ll}
$\ZZ$ & Group of relative integers \\
$\RR$ & Field of reals \\
$g$ & Integer bigger than $1$ \\
$\F$ & Free group on $2g$ generators $x_1,\ldots,x_{2g}$ \\
$R$ & Element in $\F$ given by $\prod_{j=1}^g\lbrack x_{2j-1},x_{2j}\rbrack$ \\
$\mathbf{\Pi}$ & Quotient of $\F$ by the relation $\prod_{j=1}^g\lbrack x_{2j-1},x_{2j}\rbrack$=1 \\
$\Sigma$ & Closed Riemann surface of genus $g$ \\
$\Sigma_0$ & $\Sigma$ with the interior of a disc removed \\
$K$ & Simply connected compact Lie group \\
$Z(K)$ & Centre of $K$ \\
$K_c$ & Quotient of $K$ by its centre $Z(K)$
\end{tabular}

\section*{Acknowledgement}

I would like to thank Frances Kirwan for her comments on a preliminary version of this article.
I also thank Lisa Jeffrey for sending me a copy of Shulman~\cite{S}.

The author was supported by a Marie-Curie Fellowship, EC Contract
Ref: HPMF-CT-2002-01850.
\section{Introduction}\label{sec:introduction-1}
All cohomology will have coefficients in $\RR$.

Let $\Sigma$ be a closed compact Riemann surface of genus $g$. Let $\Sigma_0$ be a compact
Riemann surface with boundary obtained from $\Sigma$ by removing a small disc.

Let $K$ be a compact simply connected Lie group. I denote by $K_c$ the quotient of $K$ by its
centre $Z(K)$. In other words, $K_c$ is the projectivisation of $K$. Let $\beta$ be in the centre
of $K$. Let $K$ act on itself by conjugation and on $K^{2g}$ by diagonal conjugation. The point
$\beta$ is fixed by $K$. Let $EK\lra BK$ be a classifying bundle for $K$.
If $M$ is a topological space acted on by $K$, the equivariant cohomology $H^*_K(M)$ of $M$ is
the singular cohomology of $M_K=M\times_K EK$. If $M$ is a manifold and it is acted
on smoothly by $K$, then $H^*_K(M)$ is also the cohomology of the Cartan-de Rham
complex $\Omega^*_K(M)$.

Let $\Phi$ be the equivariant map
$$\ba{ccc}
K^{2g} & \lra & K \\
(X_1,\ldots,X_{2g}) & \lmt & \prod_{i=1}^g\lbrack X_{2i-1},X_{2i}\rbrack.
\ea$$
Let $Y_\beta=\Phi^{-1}(\beta)$.
The space $\Phi^{-1}(\beta)/K$ is homeomorphic to the moduli space of flat connections on the trivial
bundle $K\times \Sigma_0\lra\Sigma_0$ (all $K$ principal bundles over $\Sigma_0$ are trivial) with
prescribed holonomy $\beta$ around the boundary $S^1$ of $\Sigma_0$.

Call $\kappa$ the restriction map
$$\kappa:  H^*_K(K^{2g})  \lra  H^*_K(Y_\beta).$$
Let $b_1,\ldots,b_r$ be primitive elements which generate $H^*(K)$ and let $c_1,\ldots,c_r$ be
their respective transgressions in $H^*(BK)$. One can extend $b_1,\ldots,b_r$ to equivariant classes
$b_1^K,\ldots,b_r^K$ on $K$. Let $pr_i:K^{2g}\lra K$ be the projection on the $i$-th factor and
$b_{j,i}^K$ the pull-back of $b_j^K$ by $pr_i$.

Because $\Phi$ is a moment map on the quasi-Hamiltonian space $K^{2g}$ (see \cite{AMM}), the paper
\cite{BTW} gives a way of constructing equivariant classes $a_1,\ldots,a_r$ on $Y_\beta$ (see
Equation~(\ref{eqn:definition-of-the-aj}) on page \pageref{eqn:definition-of-the-aj}) such
that $\{\kappa(c_j),\kappa(b_{j,i}^K),a_j,\\ \,j=1,\ldots,r,\,
i=1,\ldots,2g\}$ generate the $K$-equivariant cohomology of $Y_\beta$.

Also, one can construct
a $K$-equivariant $K_c$-principal bundle over $\Sigma\times Y_\beta$: the so called universal bundle. The
K\"unneth decomposition of the equivariant characteristic classes of this bundle allows me to construct
classes $\{c_j^\prime,b_{j,i}^\prime,a_j^\prime,\,j=1,\ldots,r,\,
i=1,\ldots,2g\}$ in $H^*_K(Y_\beta)$ (see Equation~(\ref{eq:definition-of-the-ajprime-bjprime-cjprime}) on
page \pageref{eq:definition-of-the-ajprime-bjprime-cjprime}). At least when $K$ is the special
unitary group and $\beta$ is a generator
of the centre of $K$, these classes are known to generate $H^*_K(Y_\beta)$ as a ring. They are the
Atiyah-Bott-Biswas-Raghavendra classes\footnote{Atiyah and Bott \cite{AB} were the first ones to construct generators
of the cohomology of the moduli space. Their classes depended on the choice of a universal bundle which was
only defined up to the tensor product with a line bundle. The universal bundle of this article is actually
the projectivisation
of theirs. It is well-defined and so are the generators constructed from it. This was proved by Biswas and
Raghavendra in \cite{BR}}.

The aim of this paper is to prove the following Theorem.
\vspace*{.2cm}
\\
\textbf{\textit{Theorem 3.2.\,---\,}}\it
Let $K$ be a compact, connected and simply connected Lie group. Let $\beta$ be in the centre of $K$ such
that $K$ acts locally freely on $Y_\beta$. Let $a_j^\prime,b_{j,i}^\prime,
c_j^\prime,a_j,b_{j,i}^K,c_j$ be equivariant classes defined as above.
The following relations hold
$$\ba{ccc}
c_j^\prime  & = & \kappa(c_j) \\
b_{j,i}^\prime  & = & \kappa(b_{j,i}^K).
\ea$$
Also, the classes $\{a_j\}$ depend on certain choices but one can choose them so that
$$a'_j=a_j.$$
\rm

This article is organised as followed. In Section~\ref{sec:gener-equiv-cohom}, I recall what
\cite{BTW} says about the restriction map involved in the reduction of $K^{2g}$ at $\beta$. I also show how to
use this result to construct generators for the $K$-equivariant cohomology of $Y_\beta$. In
Section~\ref{sec:some-class-class}, I describe the construction of equivariant characteristic
classes on $Y_\beta$ using a universal bundle over $\Sigma\times Y_\beta$. I also outline the proof of
Theorem~\ref{my_theorem}, leaving the most technical details for
Section~\ref{sec:proof-lemma}.

\section{Generators for the equivariant cohomology of $Y_\beta$}\label{sec:gener-equiv-cohom}
Let $\mathbf{\Pi}$ be the fundamental group of $\Sigma$. Fix a presentation
$$\mathbf{\Pi}=\langle x_1,\ldots,x_{2g};\, \prod_{j=1}^{g}[x_{2j-1},x_{2j}\rbrack =1\rangle.$$
Let $\F$ be the fundamental group of $\Sigma_0$. Fix a presentation
$$\F=\langle x_1,\ldots,x_{2g}\rangle$$
such that the map
$$\F\lra\mathbf{\Pi}$$
induced by the inclusion of $\Sigma_0$ in $\Sigma$, is given by the obvious map $x_k\lmt x_k$.

Let $\beta$ be in the centre of $K$. Let $Y_\beta$ be the subset of $K^{2g}$ defined by
$$Y_\beta=\left\{(X_1,\ldots,X_{2g})\in K^{2g}\mid\,\prod_{j=1}^{g}[X_{2j-1},X_{2j}\rbrack =\beta
\right\}.$$
The group $K$ acts on $K^{2g}$ by conjugation. This action restricts to $Y_\beta$, since
$\beta$ is in the centre of $K$. The quotient $\M_\beta=Y_\beta/K$ can be identified with a
moduli space of flat connections on the trivial principal bundle
$K\times\Sigma_0\lra\Sigma_0$ (because $K$ is simply connected, all principal bundles are
trivial over $\Sigma_0$), with holonomy $\beta$ around the boundary $\partial\Sigma_0$. If $K$ acts
locally freely on $Y_\beta$, then the cohomology of $\M_\beta$ is isomorphic to the
$K$-equivariant cohomology of $Y_\beta$ (recall that cohomology is taken with coefficients in $\RR$).

Let $\kappa$ be the restriction map
$$\kappa:H^*_K(K^{2g})\lra H^*_K(Y_\beta).$$
The space $K^{2g}$ has the structure of a quasi-Hamiltonian space (see \cite{AMM}) with moment map
$$\ba{cccc}
\phi: & K^{2g} & \lra & K \\
 & (X_1\ldots,X_{2g}) & \lmt & \prod_{j=1}^{g}[X_{2j-1},X_{2j}\rbrack .
\ea$$
The map $\phi$ is a submersion at any point of $Y_\beta=\phi^{-1}(\beta)$ if and only if $K$ acts
locally freely on $Y_\beta$. In this nice case, the space $\M_\beta$ is a compact symplectic orbifold
since it is the regular reduction of $K^{2g}$ at $\beta$.

Let $b_1,\ldots,b_r$ be primitive elements in $H^*(K)$ that generate the cohomology of $K$ as a
ring
$$H^*(K)=\bigwedge(\sum_{j=1}^r\RR b_j).$$
Each $b_j$ is of odd degree.

Let $c_1,\ldots,c_r$ be the transgressions in $H^*(BK)$ of respectively $b_1,\ldots,b_r$, thus
$\deg\, c_j=\deg\, b_j+1$ and
$$H^*(BK)=\RR\lbrack c_1,\ldots,c_r\rbrack .$$
For each $K$-principal bundle $G\lra X$, the classes $c_1,\ldots,c_r$ define characteristic
classes $c_1(G),\ldots,c_r(G)$ in $H^*(X)$.

\begin{prop}\label{triviality_of_the_cohomology}
Let $EK\lra BK$ be the classifying bundle for $K$. The bundle $EK\times_K K^n\lra BK$ is cohomologically
trivial as a ring, that is
$$H^*_K(K^n)\simeq H^*(BK)\otimes H^*(K^n)$$
as rings.
\end{prop}
This result is not new. I propose here a proof which is original as far as I know.

Also, because I want to state the principal results as soon as possible, I will make use in the following proof,
as well as in others, of Lemma~\ref{lemma_cj_D} even if I haven't proved it yet.
\begin{proof}
I will prove the case $n=1$. The general case is similar, one only has to use
Proposition~\ref{prop:characteristic_class_of_F} instead of Lemma~\ref{lemma_cj_D} (or use the
fact that if $X$ has an equivariantly formal action of $K$ then the diagonal action of $K$ on $X^n$ is
also equivariantly formal).

Let me first prove that $H^*_K(K)$ and $H^*(BK)\otimes H^*(K)$ are isomorphic as $H^*(BK)$ modules.
By the Leray-Hirsch Theorem,
all I have to do is prove that the homomorphism
$$H^*_K(K)\lra H^*(K)$$
is surjective. Let $b_j$ be one of the generators of $H^*(K)$. Let $D\lra S^1\times K$ be the 
$K$-equivariant $K$-principal bundle
of Lemma~\ref{lemma_cj_D}. As a $K$-principal bundle, its characteristic class $c_j(D)$ is
$\d t\otimes b_j$ (Lemma~\ref{lemma_cj_D}). As an immediate corollary, I deduce that there exists
$b_j^K$ in $H^*_K(K)$ such that $c_j(D\times_K EK)=1\otimes c_j +\d t\otimes b_j^K$, where
$b_j^K$ restricts to $b_j$. Now, because the $b_j$'s generate $H^*(K)$ as a ring, the homomorphism
$$H^*_K(K)\lra H^*(K)$$
is surjective and I have proved the existence of a $H^*(BK)$ module isomorphism between $H^*_K(K)$
and $H^*(BK)\otimes H^*(K)$. The existence of
a ring isomorphism follows because the cohomology $H^*(K)$ of the fibre in $K\times_K EK\lra BK$
is generated by classes of odd degree (see~\cite{R}).
\end{proof}

\begin{rema}\label{rema:the-section-s-and-the-pullback-of-bjK}
Let $s$ be the section
$$\ba{cccc}
s: & BK & \lra & EK\times_K K \\
 & \lbrack u\rbrack & \lmt & \lbrack u,1\rbrack .
\ea$$
Via $\id\times s:S^1\times BK\lra S^1\times K_K$, the bundle $D_K\lra S^1\times K_K$ pulls back to
$S^1\times EK\lra S^1\times BK$. Hence $s^*b_j^K=0$.
\end{rema}

Let $n$ be an integer (I will need it to be either $1$ or $2g$). Let $b_{j,i}^K$ be the pull-back,
under the projection on the $i$-th factor $K^n\lra K$, of $b_j^K$.

The following Lemma is a consequence of Proposition~\ref{triviality_of_the_cohomology}.
\begin{lemm}
The $K$-equivariant cohomology of $K^n$ is isomorphic as a ring to
$$H^*_K(K^n)\simeq \bigwedge(\sum\RR b_{j,i}^K)\otimes \RR\lbrack c_1,\ldots,c_r\rbrack .$$
\end{lemm}

Before stating the result of \BTW for $\kappa$, I need to construct some classes in $H^*_K(Y_\beta)$.

\begin{lemm}\label{the_pull-back_of_bj_is_zero}
For all $j$, the pull-back under $\phi^*$ of $b_j^K$ is zero
$$\phi^*(b_j^K)=0.$$
\end{lemm}
In their article \cite{BTW}, the authors already gave this result. I give here a proof using
Lemma~\ref{lemma_cj_D}.
\begin{proof}
Consider the diagram
$$\xymatrix{
S^1\times K^{2g} \ar^{m_1}[r]\ar^{m_2}[d] & S^1\times K \\
\Sigma_0\times K^{2g},
}$$
where the horizontal map $m_1$ is induced by the moment map $\phi:K^{2g}\lra K$ and the vertical
map $m_2$ is given by the injection of $S^1$ as the boundary of $\Sigma_0$.
The pull-back of $D$ by $m_1$ and the pull-back of
$F$ by $m_2$ are isomorphic. Indeed they are both equivariantly isomorphic to the bundle
$$(\RR\times K^{2g}\times K)/\ZZ \lra S^1\times K^{2g}$$
where the group $\ZZ$ acts by
$$\ba{ccc}
\ZZ\times\RR\times K^{2g}\times K & \lra & \RR\times K^{2g}\times K \\
(\lambda,t,(X_1,\ldots,X_{2g}),k) & \lmt & (t+\lambda,(X_1,\ldots,X_{2g}),(\prod_{j=1}^g\lbrack
X_{2j-1},X_{2j}\rbrack)^\lambda k).
\ea$$
One then deduces the Lemma from Proposition~\ref{prop:characteristic_class_of_F} and
Lemma~\ref{lemma_cj_D}
by noticing that the injection of $S^1$ as the boundary of $\Sigma_0$ induces
the null map in degree $1$ cohomology.
\end{proof}

Let $\bar{b_1},\ldots,\bar{b_r}$ be equivariant forms representing respectively $b_1^K,\ldots,b_r^K$.
\begin{rema}\label{rema:the-pull-back-of-bjK-by-s-vanish}
Such
forms will be chosen in Section~\ref{sec:proof-lemma} so that their pull-backs
by the section $s$ of Remark~\ref{rema:the-section-s-and-the-pullback-of-bjK} vanish.
\end{rema}
The
preceding Lemma implies there exist equivariant forms $\bar{a_1},\ldots,\bar{a_r}$ on $K^{2g}$ such that
\begin{eqnarray}
\phi^*(\bar{b_j}) & = & \d \bar{a_j},\,\forall j.\label{eqn:definition-of-the-aj}
\end{eqnarray}
The $\bar{a_1},\ldots,\bar{a_r}$ are closed when restricted to $Y_\beta$. They define classes
$a_1,\ldots,a_r$ in $H^*_K(Y_\beta)$.

\begin{theo}[\BTW]\label{BTW_theorem}
The equivariant cohomology $H^*_K(Y_\beta)$ is generated as a ring by the
image of $\kappa$ and the $a_1,\ldots,a_r$.
\end{theo}

The following Corollary is an immediate consequence.
\begin{coro}\label{generators_for_the_cohomology_of_Ybeta}
The equivariant cohomology $H^*_K(Y_\beta)$ is generated as a ring by the classes
$\kappa(c_j),\kappa(b_{j,i}^K),a_j$ for $j=1,\ldots,r$ and $i=1,\ldots,2g$.
\end{coro}

\section{Some classical classes in $H^*_K(Y_\beta)$}\label{sec:some-class-class}

Let $\widetilde{\Sigma}$ be the universal cover of $\Sigma$. The group $\mathbf{\Pi}$ acts on $\widetilde{\Sigma}$.
Since I have chosen a presentation for $\F$, there is a preferred isomorphism
$$\ba{ccc}
\Hom(\F,K) & \lra & K^{2g} \\
\rho & \lmt & (\rho(x_1),\ldots,\rho(x_{2g})).
\ea$$
I will use this isomorphism to identify $\Hom(\F,K)$ and $K^{2g}$. Define an action of $\mathbf{\Pi}$
on $\widetilde{\Sigma}\times Y_\beta\times K_c$
$$\ba{ccc}
\mathbf{\Pi}\times(\widetilde{\Sigma}\times Y_\beta\times K_c) & \lra & \widetilde{\Sigma}\times Y_\beta\times K_c \\
(p, \sigma, \rho, k) & \lmt & (p\sigma, \rho, \rho(\tilde{p})k),
\ea$$
where $\tilde{p}$ is any element in $\F$ lying above $p$. This action is well-defined and free. The quotient $E$,
together with the projection $\pi:E\lra\Sigma\times Y_\beta$ define a $K_c$-principal bundle. The actions of
$K$ on $Y_\beta$ by conjugation and on $K_c$ by multiplication on the left, make $E$ into a $K$-equivariant
$K_c$-principal bundle over $\Sigma\times Y_\beta$. I call it the universal bundle (this bundle is also constructed
in Jeffrey~\cite{Jef}).

\begin{exam}\label{the_case_of_SUN}
For $K=\Sun$ and $\beta$ a generator of $Z(\Sun)$, this is the projectivisation of a universal bundle for
the moduli space of semi-stable vector bundles of fixed rank and determinant (see \cite{AB} and \cite{R}).
\end{exam}

Because $K$ is compact and simply connected, I know that
$$H^*(BK)\simeq H^*(BK_c).$$
The classes $c_1,\ldots,c_r$ in $H^*(BK_c)$ define equivariant characteristic classes for the bundle
$E$ (recall that these are the same as the usual characteristic classes of the $K$-principal bundle
$E\times_K EK\lra \Sigma\times Y_\beta\times_K EK$).
These classes live in $H^*_K(\Sigma\times Y_\beta)\simeq H^*(\Sigma)\otimes H^*_K(Y_\beta)$.
Let $\alpha_1,\ldots,\alpha_{2g}$ be a basis of $H^*(\Sigma)$, dual to $x_1,\ldots,x_{2g}$ in
$\mathbf{\Pi}=\pi_1(\Sigma)$. Let $\vol$ be a volume form of volume $1$ on $\Sigma$.
Using the K\"unneth decomposition, write
\begin{eqnarray}
c_j(E\times_K EK) & = & 1\otimes c_j^\prime+\sum_{i=1}^{2g}\alpha_i\otimes b_{j,i}^\prime+\vol\otimes a_j'.\label{eq:definition-of-the-ajprime-bjprime-cjprime}
\end{eqnarray}
In Example \ref{the_case_of_SUN}, these classes are the known \ABBR generators of $H^*_K(Y_\beta)$.
A natural question is: what is the relation between these classes and the ones of
Corollary~\ref{generators_for_the_cohomology_of_Ybeta}? The answer is given in the following Theorem.

\begin{theo}\label{my_theorem}
Let $K$ be compact, connected and simply connected Lie group. Let $\beta$ be in the centre of $K$. Assume
that $\beta$ is a regular value of the moment map $\Phi$. Let the $\{a_j^\prime,b_{j,i}^\prime,
c_j^\prime,a_j,b_{j,i}^K,c_j\}$ be equivariant classes defined as above. The following relations hold
$$\ba{ccc}
c_j^\prime  & = & \kappa(c_j) \\
b_{j,i}^\prime  & = & \kappa(b_{j,i}^K).
\ea$$
Also, the classes $\{a_j\}$ depend on certain choices but one can choose them so that
$$a'_j=a_j.$$
\end{theo}
\begin{proof}
Fix $j$. I will start by proving the first part of the Theorem, that is I will
compute the $\{c^\prime_j,b^\prime_{j,i}\}$.

Let $\widetilde{\Sigma}_0$ be the universal covering for $\Sigma_0$. Define an action of
$\F=\pi_1(\Sigma_0)$ on $\widetilde{\Sigma}_0\times K^{2g}\times K$
$$\ba{ccc}
\F\times \widetilde{\Sigma}_0\times K^{2g}\times K & \lra & \widetilde{\Sigma}_0\times K^{2g}\times K \\
(p,\sigma,\rho,k) & \lmt & (p\sigma, \rho,\rho(p)k).
\ea$$
This action is free and the quotient $F$ is a $K$-principal bundle on $\Sigma_0\times K^{2g}$.
The action of $K$ by conjugation on $K^{2g}$ and by multiplication on the left on $K$ makes $F$
into a $K$-equivariant $K$-principal bundle. The situation here is somehow confusing since
there are two different actions of $K$ on $F$. In order to clarify the situation, I will call
the action that makes $F$ a $K$-principal bundle the principal action and call the action
that makes it $K$-equivariant the covering action (`covering' because it covers the action of $K$ on $K^{2g}$).

\begin{prop}\label{prop:relation_between_E_and_F}
The projectivisation of $F\mid_{\Sigma_0\times Y_\beta}$ is $K$-equivariantly isomorphic to
$E\mid_{\Sigma_0\times Y_\beta}$.
\end{prop}
\begin{proof}
The proof is easy.
\end{proof}

Call $D_n$ the restriction of $F$ to the $n$-th circle in the wedge product times
the $n$-th copy of $K$ in $K^{2g}$. I identify the circle $S^1$ with
the quotient $\RR/\ZZ$. I denote by $t$ the natural coordinate in $\RR$, thus $dt$ defines a volume form on $S^1$.

\begin{lemm}\label{lemma_Fn}
All the $D_n$'s are isomorphic to the same bundle $D\lra S^1\times K$. The total space $D$ is
$(\RR\times K\times K)/\ZZ$ where the action of $\ZZ$ is defined by
$$\ba{ccc}
\ZZ\times (\RR\times K\times K) & \lra & \RR\times K\times K \\
(\lambda,t,k_1,k_2) & \lmt & (t+\lambda, k_1, k_1^\lambda k_2)
\ea$$
and the projection is
$$\ba{ccc}
(\RR\times K\times K)/\ZZ & \lra & S^1\times K \\
\{t,k_1,k_2\} & \lmt & (\{t\},k_1).
\ea$$
\end{lemm}
\begin{proof}
Fix an integer $n$ between $1$ and $2g$ and let $x=x_n$ be the $n$-th generator of $\F$.\\
Recall that $F$ was defined as a quotient
$$F=(\widetilde{\Sigma}_0\times K^{2g}\times K)/\pi_1(\Sigma_0).$$
Let $h$ be the group homomorphism
$$\ba{ccc}
\ZZ & \lra & \pi_1(\Sigma_0) \\
\lambda & \lmt & x^\lambda.
\ea$$
This morphism makes $\ZZ$ into a subgroup of $\pi_1(\Sigma_0)$. Choose a point $\widetilde{\sigma_0}$ in
$\Sigma_0$ lying above the base point $\sigma_0$ of $\Sigma_0$. Choose a loop $S^1\lra \Sigma_0$
representing $x$. This loop lifts to a map $\chi:\RR \lra \widetilde{\Sigma}_0$ such that
$\chi(0)=\widetilde{\sigma_0}$ and for $\lambda$ an integer $\chi(t+\lambda)=x^\lambda\cdot\chi(t)$.

The covering action of $K$ (the one that covers the action of $K$ on $S^1\times K$) on $D_n$ is
induced by the conjugation on the first factor $K$ in $\RR\times K\times K$
and by multiplication on the left on the second factor $K$.

Let $\iota:K\lra K^{2g}$ be the injection of the $n$-th factor. Define a map
$$\ba{ccc}
(\RR\times K\times K)/\ZZ & \lra & (\widetilde{\Sigma}_0\times K^{2g}\times K)/\pi_1(\Sigma_0) \\
\lbrack t,X,k\rbrack & \lmt & \lbrack\chi(t), \iota(X), k\rbrack.
\ea$$
This map is well-defined since it sends the translated $\lbrack t+\lambda, X, X^\lambda k\rbrack$ of
$\lbrack t,X,k\rbrack$ by $\lambda$ to
$$\ba{cl}
 & \lbrack\chi(t+\lambda), \iota(X), X^\lambda k\rbrack \\
= & \lbrack x_n^\lambda\cdot\chi(t),\iota(X), X^\lambda k\rbrack \\
= & \lbrack\chi(t), \iota(X), k\rbrack.
\ea$$
It is also $K$-equivariant because, for $\ell$ in $K$, it sends $\lbrack t,\Ad_\ell X, \ell k\rbrack $ to
$$\ba{cl}
 & \lbrack\chi(t), \iota(Ad_\ell X), \ell k\rbrack  \\
= & \lbrack\chi(t), Ad_\ell\circ \iota(X), \ell k\rbrack  \\
= & \ell\cdot\lbrack\chi(t), \iota(X), k\rbrack .
\ea$$
The Lemma follows.
\end{proof}

\begin{lemm}\label{lemma_cj_D}
The characteristic classes of $D$ are
$$c_j(D)=dt\otimes b_j.$$
\end{lemm}
\begin{proof}
The proof is to be found in the next Section.
\end{proof}

The previous Lemma is actually very useful since it can be used to prove Proposition
\ref{prop:characteristic_class_of_F} which leads to Theorem \ref{my_theorem}, and also
Proposition \ref{triviality_of_the_cohomology} and Lemma
\ref{the_pull-back_of_bj_is_zero}.

Now I wish to compute the equivariant characteristic classes of $F$.
The cohomologies of $\Sigma$ and $\Sigma_0$ are the same in degree $0$ and $1$ and $\Sigma_0$ has no
cohomology in degree $2$. I use the same basis $\alpha_1,\ldots,\alpha_{2g}$ for $H^1(\Sigma)$ and
$H^1(\Sigma_0)$.

\begin{prop}\label{prop:characteristic_class_of_F}
Using the K\"unneth decomposition, the equivariant characteristic classes of $F$ are
$$c_j(F\times_K EK)=1\otimes c_j + \sum_{i=1}^{2g}\alpha_i\otimes b_{j,i}^K.$$
\end{prop}

\begin{proof}
The surface $\Sigma_0$ is homotopic to a wedge product of $2g$ circles. As is easily seen,
Proposition \ref{prop:characteristic_class_of_F} will follow from the computation of the characteristic classes
of the restriction of $F$ to the $n$-th circle in the wedge product times the $m$-th copy of $K$
in $K^{2g}$, that is $D_n$. When $n\neq m$, the restriction of $F$ is trivial.
Thus, I only have to do the computation for $n=m$.
This was done in Lemma~\ref{lemma_cj_D}.
\end{proof}

The first part of the Theorem is a consequence of Proposition~\ref{prop:relation_between_E_and_F} and
Proposition~\ref{prop:characteristic_class_of_F}. Indeed, by Proposition~\ref{prop:relation_between_E_and_F}
$$c_j(E\times_K EK)\mid_{\Sigma_0\times Y_\beta}=c_j(F\times_K EK)\mid_{\Sigma_0\times Y_\beta}.$$
But
\begin{eqnarray}
c_j(E\times_K EK)\mid_{\Sigma_0\times Y_\beta}=
1\otimes c_j^\prime+\sum_{i=1}^{2g}\alpha_i\otimes b_{j,i}^\prime \label{premiere_equation}
\end{eqnarray}
and by Proposition~\ref{prop:characteristic_class_of_F}
\begin{eqnarray}
c_j(F\times_K EK)\mid{\Sigma_0\times Y_\beta} & = & (\id^*\times\kappa)c_j(F\times_K EK) \nonumber\\
 & = & 1\otimes \kappa(c_j)+\sum_{i=1}^{2g}\alpha_i\otimes \kappa(b_{j,i}^K). \label{eqn:deuxieme_equation}
\end{eqnarray}
The first part of the Theorem follows by comparing Equation~(\ref{premiere_equation})
and Equation~(\ref{eqn:deuxieme_equation}).

The second part of Theorem~\ref{my_theorem} is a consequence of the following Lemma and of the definition
of the classes $a_1,\ldots,a_r$.
\begin{lemm}\label{lemm:extension-of-a-prime-j-to-K2g}
One can choose representatives $\bar{b}_j^K$ and $\bar{a}_j^\prime$ of respectively $b_j^K$ and
$a_j^\prime$ such that $\bar{a}_j^\prime$ extends to an equivariant class on $K^{2g}$ whose differential
is $\Phi^*(\bar{b}_j^K)$.
\end{lemm}
The proof of this Lemma, together with the proof of Lemma~\ref{lemma_cj_D}, is given in the next Section.

This finishes the proof of Theorem~\ref{my_theorem}.
\end{proof}

\begin{rema}
If the group $K$ is a torus $T$, then it is not simply connected and one can not apply
Theorem~\ref{my_theorem} to find generators of the cohomology of the moduli space. Nevertheless,
in this case the moduli space is either $T^{2g}$ if $\beta=1$ or is empty otherwise. Its
cohomology is thus very simple. In general, a connected compact Lie group is the quotient by a finite group
of a product of a simply connected group and a torus. Combining what I have just said about
tori with Theorem~\ref{my_theorem}, one can deduce generators for the cohomology of the moduli
space.
\end{rema}

\section{Proof of Lemma~\ref{lemm:extension-of-a-prime-j-to-K2g} and Lemma \ref{lemma_cj_D}}\label{sec:proof-lemma}
The aim of this section is to prove Lemma~\ref{lemma_cj_D} and Lemma~\ref{lemm:extension-of-a-prime-j-to-K2g}.

The proof of Lemma~\ref{lemma_cj_D} goes as follows: in a similar
way to Jeffrey\cite{Jef}, I construct explicitly a map from $S^1\times K$ to the fat realisation of
$BK$ as a simplicial manifold such that the pull-back of $EK$ under this map is the principal bundle $D_n$.
I then use a result of Shulman to prove Lemma~\ref{lemma_cj_D}.

I will start by introducing the simplicial realisations of $EG$ and $BG$ for any group $G$ (see~\cite{Se}). Let
$$\Delta^m=\left\{(t_0,\ldots,t_m)\in\lbrack0,1\rbrack ^{m+1}\mid\, \sum_{i=0}^m t_i=1\right\}$$
be the standard $m$-simplex. The simplicial realisation $\NGb$, respectively $\NG $, of $EG$,
respectively $BG$,  is given by
the contravariant functor
$$\{\Delta^m,\, m\in\NN\}\lra \{\mbox{category of smooth manifolds}\}$$
which to each $\Delta_m$ associates
$\NGb(m)=G^{m+1}$, respectively $\NG (m)=G^m$. For $i$ between $0$ and $m$, the map
$\bar{\epsilon_i}:\NGb(m)\lra\NGb(m-1)$,
respectively $\epsilon_i:\NG (m)\lra\NG (m-1)$, corresponding
to the $i$-th face map $\epsilon^i:\Delta^{m-1}\lra\Delta^m$ is given by the omission of the
$i$-th term
$$\bar{\epsilon_i}(k_0,\ldots,k_m)=(k_0,\ldots,\check{k}_i,\ldots,k_m),$$
respectively
\begin{eqnarray}
\quad\quad\epsilon_i(k_1,\ldots,k_m) & = & \left\{
\ba{ll}
(k_2,\ldots,k_m) & i=0 \\
(k_1,k_2,\ldots,k_i k_{i+1},\ldots,k_m) & i=1,\ldots,m-1 \\
(k_1,\ldots,k_{m-1}) & i=m.
\ea\right.
\end{eqnarray}\label{definition-of-varepsilon-i-for-NK}
The projection $\NGb(m)\lra \NG (m)$ is the map
$$q_m(k_0,\ldots,k_m)=(k_1 k_2^{-1},\ldots,k_{m-1}k_m^{-1}).$$
I then take $EG$, respectively $BG$, to be the fat realisation of $\NGb$, resp. $\NG $, that
is $EG$ is the quotient of $\bigsqcup_m \Delta^m\times\NGb(m)$ by the relation
$$(\epsilon^i\times \id)(t,x)\sim(id\times \bar{\epsilon_i})(tx)$$
for $(t,x)$ in $\Delta^{m-1}\times\NGb(m)$ and $BG$ is the quotient of
$\bigsqcup_m \Delta^m\times \NG (m)$ by the relation
$$(\epsilon^i\times \id)(t,x)\sim(id\times \epsilon_i)(tx)$$
for $(t,x)$ in $\Delta^{m-1}\times \NG (m)$. The projection $q:EG\lra BG$ is induced by the maps
$q_m$. The action of $G$ on $EG$ is induced by the diagonal multiplication of $G$ on the right
on $\NGb(m)=G^{m+1}$ for each $m$. In fact, there is a second action of $G$ on $EG$, the covering
action, which descends to an action of $G$ on $BG$. This action on $EG$ is induced by the diagonal
multiplication on the left of $G$ on $\NGb(m)=G^{m+1}$, it descends to an action on $BG$ induced
by the diagonal conjugation on $\NG (m)=G^m$.

Notice that
$$\ba{rcl}
S^1 & = & B\ZZ \\
\Sigma & = & B\mathbf{\Pi} \\
\Sigma_0 & = & B\F.
\ea$$
In the simplicial models, the injection of $\Sigma_0$ in $\Sigma$ is given by the projections
$$\F^n\lra\mathbf{\Pi}^n.$$
The surface $\Sigma_0$ is homotopic to the wedge product of $2g$ circles. The inclusion of the
$i$-th circle $S^1\lra\Sigma_0$ is given by
$$\ba{ccc}
\ZZ^n & \lra & \F^n \\
\lambda & \lmt & (x_i^\lambda,\ldots,x_i^\lambda).
\ea$$
Let $R=\prod_{j=1}^g\lbrack x_{2j-1},x_{2j}\rbrack$ be in $\F$.
The inclusion of $S^1$ as the boundary of $\Sigma_0$ is the realisation of the simplicial map
$$\ba{ccc}
\ZZ^n & \lra & \F^n \\
\lambda & \lmt & (R^\lambda,\ldots,R^\lambda).
\ea$$
The classifying map $f:S^1\times K\lra BK$ for the bundle $D$ is the realisation of the simplicial map
$$\ba{cccc}
f: & \ZZ^n\times K & \lra & K^n \\
 & (k,\lambda_1,\ldots,\lambda_n) & \lmt & (k^{\lambda_1},\ldots,k^{\lambda_n}).
\ea$$
The homology of $S^1$ is the Eilenberg-Mac Lane group homology $H_*(\ZZ)$ of $\ZZ$ (see \cite{Mac}). Denote by
$x$ a generator of $\ZZ$. A generator of $H_1(\ZZ)$ is given by $x$.

Since the centre of $K$ is finite, on the one hand $K$ and $K_c$, and on the other hand $BK$ and $BK_c$ have naturally
isomorphic cohomologies. In addition, every
$K$-principal bundle defines a $K_c$-principal bundle in a natural way: one just has to compose the classifying
map of the $K$-principal bundle with the projection $BK\lra BK_c$. The simplicial version of this last projection
is given by the maps $K^n\lra K_c^n$. Moreover, the characteristic classes of the induced $K_c$-principal bundle
are the same as the one of the $K$-principal bundle.

The cohomology of $BK_c$ is the cohomology of a
double complex $\Omega^{p,q}=\Omega^p(\NK_c(q))$
with two differentials: one is the usual differential of
forms $\d:\Omega^p(\NK_c(q))\lra \Omega^{p+1}(\NK_c(q))$
and the other is
$$\ba{c}
\delta:  \Omega^p(\NK_c(q))  \lra  \Omega^p(\NK_c(q+1))\\
\delta=\sum_{i=0}^{q+1}(-1)^i\varepsilon_i^*
\ea$$
Recall that the maps $\varepsilon_i$ were defined in~(\ref{definition-of-varepsilon-i-for-NK}). The two
differentials $\d$ and $\delta$ commute and the total differential is $\d+(-1)^p\delta$ on $\Omega^{p,q}$.
Shulman~\cite{S} (see also Bott, Shulman and Stasheff \cite{BSS}) proved that each class $c_j$
in $H^{2r_j}(BK_c)$ can be represented by a form
$c_{j,1}+\ldots+c_{j,r}$ in $\Omega^{*,*}$, with $c_{j,k}\in \Omega^{2r-k,k}$. That this form is closed
means that
$$\ba{rcl}
\d c_{j,1} & = & 0 \\
\delta c_{j,k} & = & (-1)^k\d c_{j,k+1}\mbox{ for }1\leq j<r_j \\
\delta c_{j,r} & = & 0.
\ea$$
Also, the closed form $c_{j,1}\in \Omega^{2r_j-1}(K_c)$ is a representative of $b_j$.

Everything in the above paragraph can be done equivariantly. One just has to replace cohomology by
equivariant cohomology, the complex $\Omega^{p,q}=\Omega^p(\NK_c(q))$ by
$\Omega^{p,q}_K=\Omega^p_K(\NK_c(q))$,
the differential $\d$ by $\d_K$, the forms $c_j$ by their respective extensions $c_j^K$ in $H_{K}^{2r_j}(BK_c)$,
and so on. Hence it is clear
that the pairing $\langle x,f^* c_j\rangle$ is $b^K_j$. This together with the fact that the restriction of
$D$ to $S^1$ is trivial proves Lemma~\ref{lemma_cj_D}.

I chose $c_{j,1}^K$ as a representative $\bar{b}_j^K$ of $b_j^K$. The section
$$\ba{cccc}
s: & BK & \lra & EK\times_K K \\
 & \lbrack u\rbrack & \lmt & \lbrack u,1\rbrack
\ea$$
is obtained from the inclusion of the identity $1$ into $K$ by applying the functor, from $K$-spaces to spaces,
$$X\lmt X\times_K EK.$$
Because the inclusion map $\{1\}\lra K$ is constant and because $c_{j,1}^K$ is of odd degree, it follows that
$s^*b_j^K=s^*c_{j,1}^K=0$. This proves Remark~\ref{rema:the-pull-back-of-bjK-by-s-vanish}.

The equivariant forms $c_{j,1}^K$ and $c_{j,2}^K$ satisfy
$\d_K c_{j,2}^K=-\delta c_{j,1}^K$. A classifying map for the $K$-equivariant $K_c$-principal bundle
$E$ is the realisation of the equivariant simplicial map
$$\ba{cccc}
\Psi: & Y_\beta\times\Pi^n & \lra & K_c^n \\
 & (\rho,p_1,\ldots,p_n) & \lmt & (\rho(p_1),\ldots,\rho(p_n)).
\ea$$
The homology of $\Sigma$ is the Eilenberg-Mac Lane group homology $H_*(\mathbf{\Pi})$. A generator of
$H_2(\mathbf{\Pi})$ is (see \cite[Proposition 3.9]{G})
$$c=\sum_{i=1}^{2g}\frac{\del R}{\del x_i}\otimes x_i,$$
where $\frac{\del R}{\del x_i}$ is the Fox free differential calculus (see \cite[Section 3]{G}).
In fact,
$$\frac{\del R}{\del x_i}=\gamma^0_j-\gamma^1_i,$$
where $\gamma^\tau_i\in\F$ is
$$\ba{ccc}
\gamma^0_{2i-1} & = & \prod_{l=1}^{i-1}\lbrack x_{2l-1},x_{2l}\rbrack, \\
\gamma^1_{2i-1} & = & \gamma^0_{2i-1}x_{2i-1}x_{2i}x_{2i-1}^{-1}, \\
\gamma^0_{2i} & = & \gamma^0_{2i-1}x_{2i-1}, \\
\gamma^1_{2i} & = & \gamma^0_{2i-1}\lbrack x_{2i-1},x_{2i}\rbrack.
\ea$$
A representative $\bar{a}_j^\prime$ of $a_j^\prime$ is given by
$$\ba{rcl}
\bar{a}_j^\prime & = & \langle\sum_{i=1}^{2g}\frac{\del R}{\del x_i}\otimes x_i,\Psi^*c_{j,2}^K\rangle \\
 & = & \sum_{i=1}^{2g}\sum_{\tau=0,1}(-1)^\tau\langle\gamma^\tau_i\otimes x_i,\Psi^*c_{j,2}^K\rangle \\
 & = & \sum_{i=1}^{2g}\sum_{\tau=0,1}(-1)^\tau(\ev_{\gamma^\tau_i}\times\ev_{x_i})^*c_{j,2}^K,
\ea$$
where, for $p\in\F$, the evaluation map is
$$\ba{cccc}
ev_p: & \Hom(\F,K_c) & \lra & K_c \\
 & \rho & \lmt & \rho(p).
\ea$$
The last line of the above computation of $\bar{a}_j^\prime$ clearly shows that $\bar{a}_j^\prime$
extends to a form on $K^{2g}$. The differential of this form is
\begin{eqnarray}
\d\bar{a}_j^\prime & = & \sum_{i=1}^{2g}\sum_{\tau=0,1}(-1)^\tau(\ev_{\gamma^\tau_i}\times\ev_{x_i})^*\d c_{j,2}^K \nonumber\\
 & = & -\sum_{i=1}^{2g}\sum_{\tau=0,1}(-1)^\tau(\ev_{\gamma^\tau_i}\times\ev_{x_i})^*\delta c_{j,1}^K \nonumber\\
&=&-\sum_{i=1}^{2g}\sum_{\tau=0,1}(-1)^\tau(\ev_{\gamma^\tau_i}\times\ev_{x_i})^*\sum_{l=0}^2(-1)^l\varepsilon_l^*c_{j,1}^K\nonumber\\
&=&-\sum_{i=1}^{2g}\sum_{\tau=0,1}\sum_{l=0}^2(-1)^{\tau+l}(\varepsilon_l\circ(\ev_{\gamma_i^\tau}\times \ev_{x_i}))^*c_{j,1}^K \nonumber\\
&=&\sum_{i=1}^{2g}\sum_{\tau=0,1}\sum_{l=0}^2(-1)^{1+\tau+l}\ev_{z^\tau_{i,l}}^*c_{j,1}^K,\label{eqn-last-line}
\end{eqnarray}
where
$$z^\tau_{i,0}=x_i,\,z^\tau_{i,1}=\gamma^\tau_i x_i\mbox{ and } z^\tau_{i,2}=\gamma^\tau_i.$$
Because
$$\ba{rcl}
\ev^*_{z^0_{i,0}} & = & \ev^*_{z^1_{i,0}} \\
\ev^*_{z^1_{2i,1}} & = & \ev^*_{z^1_{2i-1,2}} \\
\ev^*_{z^0_{2i-1,1}} & = & \ev^*_{z^0_{2i,2}} \\
\ev^*_{z^0_{2i,1}} & = & \ev^*_{z^1_{2i-1,1}} \\
\ev^*_{z^0_{2i+1,2}} & = & \ev^*_{z^1_{2i,2}}
\ea$$
many terms cancel each other in the expression (\ref{eqn-last-line}) above. I am left with
$$\d\bar{a}_j^\prime=\ev_{z^1_{2g,2}}^*c_{j,1}^K-ev_{z^0_{1,2}}^*c_{j,1}^K.$$
But $z^0_{1,2}=1$ the identity of $\F$ and hence $ev_{z^0_{1,2}}$ is a constant map, whereas
$z^1_{2g,2}=\prod_{j=1}^g\lbrack x_{2j-1},x_{2j}\rbrack$ and hence
$\ev_{z^1_{2g,2}}^*=\Phi^*$. Since $c_{j,1}^K$ is a representative of $b_j^K$, I deduce that on $K^{2g}$
$$\d\bar{a}_j^\prime=\Phi^* \bar{b}_j^K.$$
This proves Lemma~\ref{lemm:extension-of-a-prime-j-to-K2g}.

{\providecommand{\bysame}{\leavevmode ---\ }}

\nobreak
\end{document}